\documentclass[12pt,twoside]{article}
\usepackage{amsmath,amsbsy,amsfonts}

\newtheorem{theorem}{Theorem}[section]
\newtheorem{corollary}{Corollary}[section]
\newtheorem{claim}{Claim}[section]

\let\Section=\section
\def\section{\setcounter{equation}{0}\Section}

\begin{document}

\date{}
\title{\textbf{Existence of blow-up solutions for a class of elliptic system
with convection term}}
\author{\textsf{{Claudianor O. Alves} \thanks{%
C.O. Alves was partially supported by CNPq/Brazil 303080/2009-4,
coalves@dme.ufcg.edu.br}} \\
{\small \textit{Unidade Acad\^emica de Matem\'atica}}\\
{\small \textit{Universidade Federal de Campina Grande}}\\
{\small \textit{58429-900, Campina Grande - PB - Brazil}}\\
{\small \textit{e-mail address: coalves@dme.ufcg.edu.br}} \\
\\
\vspace{1mm} \textsf{{Dragos-Patru Covei} }\\
{\small \textit{Department of Applied Mathematics}}\\
{\small \textit{The Bucharest University of Economics Study }}\\
{\small \textit{Piata Romana, 1st district, postal code: 010374, postal
office: 22, Romania}}\\
{\small \textit{e-mail address:coveidragos@yahoo.com}}}
\maketitle

\begin{abstract}
The present paper concerns with the existence of solutions for a class of
elliptic systems involving nonlinearities of the Keller-Osserman type and
combined with the convection terms. Firstly, we establish a result involving
sub and super-solution for a class of elliptic system whose nonlinearity can
depend of the gradient of the solution. This result permits to study the
existence of blow-up solution for a large class of systems.
\end{abstract}

\vspace{0.5 cm} \noindent \textbf{\footnotesize 2000 Mathematics Subject
Classifications:} {\scriptsize 35B44, 35B08, 35J05}.\newline

\section{Introduction}

In this article, we study the existence of solutions for the following class
of elliptic system with convection term 
\begin{equation*}
\left\{ 
\begin{array}{l}
\Delta u+b_{1}(x)|\nabla u|^{q_{1}}=F_{u}(x,u,v)\text{ in }\Omega , \\ 
\mbox{} \\ 
\Delta v+b_{2}(x)|\nabla v|^{q_{2}}=F_{v}(x,u,v)\text{ in }\Omega,%
\end{array}%
\right. \eqno{(S)}
\end{equation*}%
where $\Omega \subset \mathbb{R}^{N}(N\geq 1)$ is a bounded domain with
smooth boundary or $\Omega =\mathbb{R}^{N}$, $0<q_{i}\leq 2$, $b_{i}:%
\overline{\Omega} \rightarrow \mathbb{R}^{+}$ ($i=1,2$) are continuous
functions and $F:\mathbb{R}^{N}\times \mathbb{R}^{+}\times \mathbb{R}%
^{+}\rightarrow \mathbb{R}^{+}$ is a $C^{1}$ function verifying some
technical condition, which are mentioned later on.

For the case where $\Omega$ is a bounded domain, the system will be studied
under three different types of boundary conditions:

\begin{itemize}
\item \textbf{Finite Case:} Both components $(u,v)$ bounded on $\partial
\Omega $, that is, 
\begin{equation}
\left\{ 
\begin{array}{lll}
u=\alpha \,\,\,~~\mbox{on}~~\partial \Omega , &  &  \\ 
v=\beta \,\,\,~~\mbox{on}~~\partial \Omega , &  & 
\end{array}%
\right.  \tag{F}  \label{F}
\end{equation}%
with $\alpha ,\beta \in (0,+\infty )$.

\item \textbf{Infinite Case:} Both components blowing up simultaneously on $%
\partial \Omega $, that is, 
\begin{equation}
\left\{ 
\begin{array}{lll}
u=+\infty \,\,\,~~\mbox{on}~~\partial \Omega , &  &  \\ 
v=+\infty \,\,\,~~\mbox{on}~~\partial \Omega , &  & 
\end{array}%
\right.  \tag{I}  \label{I}
\end{equation}%
where $u=+\infty $ on $\partial \Omega $ and $v=+\infty $ on $\partial
\Omega $ should be understood as $u(x)\rightarrow +\infty $ and $%
v(x)\rightarrow +\infty $ as $\mbox{dist}(x,\partial \Omega )\rightarrow 0$.

\item \textbf{Semifinite Case:} One of the components bounded while the
other one blows up on $\partial \Omega $, that is, 
\begin{equation}
\left\{ 
\begin{array}{lll}
u=+\infty \,\,\,~~\mbox{on}~~\partial \Omega , &  &  \\ 
v=\beta \,\,\,\,\,\,\,\,\,\,~~\mbox{on}~~\partial \Omega , &  & 
\end{array}%
\right.   \tag{SF1}  \label{SF1}
\end{equation}

or

\begin{equation}
\left\{ 
\begin{array}{lll}
u=\alpha \,\,\,~~~~~~~~\mbox{on}~~\partial \Omega , &  &  \\ 
v=+\infty \,\,\,\,\,\,\,\,\,\,~~\mbox{on}~~\partial \Omega . &  & 
\end{array}%
\right.  \tag{SF2}  \label{SF2}
\end{equation}
\end{itemize}

A solution $(u,v)\in C^{2}(\Omega)\times C^{2}(\Omega)$ of the system $(S)$
is called a \textit{blow-up solution} if the condition $(I)$ holds and 
\textit{semifinite blow-up solution} when $(SF1)$ or $(SF2)$ holds.

For the case $\Omega=\mathbb{R}^{N}$, we consider the following class of
elliptic systems

\begin{equation}  \label{LS}
\left \{ 
\begin{array}{lll}
\Delta u + b_1(x)|\nabla u|^{q_1}= F_{u}(x,u,v)\,\,\, ~~ \mbox{in} ~~ 
\mathbb{R}^{N}, &  &  \\ 
\Delta v + b_2(x)|\nabla v|^{q_2}= F_{v}(x,u,v)\,\,\, ~~ \mbox{in} ~~ 
\mathbb{R}^{N}, &  &  \\ 
u, v>0 \,\,\,\,\,\,\,\,\,\,\,\,~~~~~~~~~~~~~ \mbox{in} ~~ \mathbb{R}^{N}. & 
& 
\end{array}
\right.  \tag{LS}
\end{equation}
\noindent Associated with this class of systems, our main result is
concerned with the existence of \textit{entire large solutions}, that is,
solutions $(u,v)$ satisfying $u(x)\rightarrow +\infty$ and $v(x)\rightarrow
+\infty$ as $|x|\rightarrow +\infty$.

The scalar case associated with system (S), namely 
\begin{equation*}
\Delta u+b(x)|\nabla u|^{q}=F_{u}(x,u)\,\text{in}\,\,\Omega ,
\end{equation*}%
has been considered by several authors. We would like to mention the papers
of Alarc\'on, Garcia-Melian \& Quass \cite{AMQ}, Bandle \& Giarrusso \cite%
{bg}, Bandle \& Marcus \cite{bandle}, Covei \cite{covei}, Filippucci, Pucci
\& Rigoli \cite{FPR}, Garcia \& Melian \cite{mel}, Garcia-Melian \& Rossi 
\cite{rosi}, Ghergu,Niculescu \& Radulescu \cite{GH} , Holanda \cite{holanda}%
, Lair \cite{lair}, Lair \& Wood \cite{aw}, Mohammed \cite{ahmed}, Keller 
\cite{keller}, Osserman \cite{osserman}, Mi \& Liu \cite{ML}, and references
therein.

For instance, Lair \cite{lair} showed the existence of solutions of the
problem%
\begin{equation*}
\Delta u=r\left( x\right) h\left( u\left( x\right) \right) \text{ for }x\in
\Omega \subseteq \mathbb{R}^{N}\text{, }N\geq 3
\end{equation*}%
where the function $h:\left[ 0,+\infty \right) \rightarrow \left[ 0,\infty
\right) $ satisfies the $\mathcal{F-}$condition : 
\begin{equation*}
\begin{array}{l}
\bullet \,\,\, h \in C^{1}\left( \left[ 0,\infty \right) \right) ,\text{ }%
h\left( 0\right) =0,\text{ }h^{\prime }\left( t\right) \geq 0\text{ }\forall
t\in \left[ 0,\infty \right) , \\ 
\mbox{} \\ 
\bullet \,\,\, h\left( t\right) >0\text{ }\forall t\in \left( 0,\infty
\right)%
\end{array}
\end{equation*}
and the well known Keller-Osserman condition ( \cite{keller}, \cite{osserman}
), that is, 
\begin{equation}
\bullet \,\,\,\, \int_{1}^{\infty }\frac{1}{H\left( t\right) ^{1/2}}%
dt<\infty ,\text{ \ } \;\;\;\;\;\;\;\;\; H\left( t\right)
:=\int_{0}^{t}h\left( s\right) ds,  \label{ko}
\end{equation}%
and the function $r:\overline{\Omega }\rightarrow \left( 0,\infty \right) $
satisfies the $\mathcal{P-}$condition:

\begin{equation*}
r\in C^{0,\alpha }\left( \overline{\Omega }\right) ,\,\,\,\mbox{if}%
\,\,\,\Omega \,\,\,\mbox{is a bounded domain}
\end{equation*}%
or 
\begin{equation*}
r\in C_{loc}^{0,\alpha }\left( \mathbb{R}^{N}\right) ,\,\,\,\mbox{if}%
\,\,\,\Omega =\mathbb{R}^{N}.
\end{equation*}

In \cite{aw}, Lair \& Wood proved the existence of non-negative solutions of
the problem%
\begin{equation*}
\Delta u+\left\vert \nabla u\right\vert ^{q}=r\left( x\right) u^{\gamma }%
\text{ in }\Omega
\end{equation*}%
with $r \in \mathcal{P}$, $q\in \left( 0,2\right] $ and $\gamma >\max
\left\{ 1,q\right\} $.

Later, Ghergu,Niculescu \& Radulescu \cite{GH} considered the equation 
\begin{equation*}
\Delta u+q\left( x\right) \left\vert \nabla u\right\vert ^{a}=r\left(
x\right) h\left( u\right) \text{ in }\Omega ,
\end{equation*}%
assuming $a\in \left( 0,2\right]$, $h\in \mathcal{F}$ and 
\begin{equation*}
H\left( t\right)/h^{2/a}\rightarrow 0 \,\,\, \mbox{as} \,\,\, t\rightarrow
\infty,
\end{equation*}
for some suitable functions $q$, $r$ of the class $\mathcal{P}$. An
important common point among the above papers is the fact that they assume
that nonlinearity is monotone.

Recently, Alves \& Holanda in \cite{alves}, combining variational method
with the existence of sub and super-solution, obtained solutions with
boundary conditions (\ref{F}), (\ref{I}) and (SF), for the system of the
form 
\begin{equation}
\Delta U=\nabla F\left( x,U\right) \text{ in }\Omega   \label{ao}
\end{equation}%
in which $U:=\left( u,v\right) ,$ 
\begin{equation*}
\nabla F:=\left\{ 
\begin{array}{l}
\left( F_{u}\left( x,u,v\right) ,F_{v}\left( x,u,v\right) \right) \text{ if }%
\Omega \subset \mathbb{R}^{N} \\ 
\mbox{} \\ 
\left( a_{1}\left( x\right) F_{u}\left( x,u,v\right) ,a_{2}\left( x\right)
F_{v}\left( x,u,v\right) \right) \text{ if }\Omega =\mathbb{R}^{N},%
\end{array}%
\right. 
\end{equation*}%
and for suitable functions which are not necessarily monotone, but for a
particular class of systems of the form $(S)$, where $b_{i}=0$ for $i=1,2$.

The motivation to study system $(S)$ comes from the study of the chemical,
physical, biological and economical phenomena, see \cite{C,D} for details .
Also such systems can model phenomena from the study of ecological
prey-predator models, in that context we refer to Leung \cite{LE,LE2}.

Throughout this article, we assume that $b_{i}\in \mathcal{P}$ and $%
F_{u},F_{v}$ are locally H\"{o}lder continuous with exponent $\alpha \in
\left( 0,1\right) $, verifying the following additional condition:\newline

\noindent There are $a_{i}$, $a_{i}^{2}\in \mathcal{P}$ ($i=1,2$) and $f_{i}$%
, $g\in \mathcal{F}$, satisfying%
\begin{equation}
F_{t}(x,t,s)\geq a_{1}(x)f_{1}(t)\text{ }\forall x\in \overline{\Omega },%
\text{ }t,s>0\text{ }  \label{1.2}
\end{equation}%
\begin{equation}
F_{s}(x,t,s)\geq a_{2}(x)f_{2}(s)\text{ }\forall x\in \overline{\Omega },%
\text{ }t,s>0\text{ }  \label{1.2.0}
\end{equation}%
and 
\begin{equation}
g(t)>\max_{i=1,2}\left\{ \frac{F_{i}\left( x,t,t\right) }{a_{i}^{2}(x)}%
\right\} \text{ }\forall x\in \overline{\Omega },\text{ }t>0.  \label{1.4}
\end{equation}%
A simple example of nonlinearity $F$ satisfying the above assumptions is%
\begin{equation*}
F\left( x,u,v\right) =c_{1}\left( x\right) u^{\rho }+c_{2}\left( x\right)
u^{\sigma }v^{\gamma }+c_{3}\left( x\right) v^{\theta }
\end{equation*}%
where $\rho $, $\theta >2$, $\sigma +\gamma >2$ with $\sigma ,\gamma <2$ and 
$c_{i}$ ($i=1,2,3$) are some suitable functions.

Not before to enumerate our results about the considered system we wish to
say that if the nonlinearities are not necessarily non-decreasing it is
known that the problem of uniqueness of solution is not so easy even we
refer to the scalar case treated in various references. But, for some
particular cases of nonlinearities we can see that many authors appeals to
the asymptotic behavior of the solution in order to prove the uniqueness of
explosive solutions for both scalar and systems cases. In this article we
will restrict our research only to the problem of existence of solutions.
The uniqueness problem becomes more delicate topic included in our future
goals.

Our first result related to the problem $(P)$ is the following: \ 

\begin{theorem}
\label{1.1}Suppose that $\Omega $ is a bounded domain, $b_{i}\in \mathcal{P}$
and (\ref{1.2})-(\ref{1.4}) hold. Then:

\begin{description}
\item[i)] Problem $(P)$ admits positive solution with the boundary condition %
\eqref{F}.

\item[ii)] Problem $(P)$ admits positive solution with the boundary
condition \eqref{I}.

\item[iii)] Problem $(P)$ admits positive solution with the boundary
condition (\ref{SF1}) or (\ref{SF2}).
\end{description}
\end{theorem}

Our next result is related to existence of entire large solution for system $%
(S)$ for the case where $\Omega =\mathbb{R}^{N}$. For expressing the next
result, we assume that functions $a_{i}^{2}$ ($i=1,2$) belongs to $\mathcal{P%
}$ and that the problem%
\begin{equation}
\begin{array}{c}
-\Delta z\left( x\right) =\overset{2}{\underset{i=1}{\sum }}a_{i}^{2}\left(
x\right) \text{ for }x\in \mathbb{R}^{N}\text{, \ }z\left( x\right)
\rightarrow 0\text{ as }\left\vert x\right\vert \rightarrow \infty%
\end{array}
\label{1.6.}
\end{equation}%
has a $C^{2}$ supersolution.

\begin{theorem}
\label{entire}Assume that (\ref{1.2})-(\ref{1.6.}) hold. Then system $(P)$
has an entire large solution.
\end{theorem}

Before to conclude this introduction, we would like to say that Theorems \ref%
{1.1} and \ref{entire} complete the study made in \cite{alves}, in the sense
that, in that paper the authors considered only the case where $b_i=0
(i=1,2).$ Moreover, we would like to detach that the authors does know any
result involving sub and supersolution that can be used for system $(S)$. To
overcome this difficulty, we prove in Section 2 a result that allows us to
apply sub and supersolution for $(S)$.

\section{An auxiliary system}

In this section, we will work with an auxiliary system associated with $(S)$%
. In what follows, fixed $R>0$, we denote by $\xi_R : [0,+\infty) \to
[0,+\infty)$ be a nondecreasing continuous functions verifying 
\begin{equation*}
\xi_R(t)=t \,\,\, \mbox{for} \,\,\, t \in [0,R] \,\,\, \mbox{and} \,\,\,
\xi(t)=R \,\,\, \mbox{for} \,\,\, t \geq R.
\end{equation*}
Using this function, we consider the system 
\begin{equation*}
\left\{ 
\begin{array}{l}
\Delta u +b_{1}(x) \xi_R(| \nabla u |^{q_{1}})=F_{u}(x,u,v) \text{ in }%
\Omega , \\ 
\mbox{} \\ 
\Delta v +b_{2}(x) \xi_R(| \nabla v |^{q_{2}})=F_{v}(x,u,v) \text{ in }%
\Omega, \\ 
\mbox{} \\ 
u=\alpha, v=\beta \,\, \mbox{on} \,\, \partial \Omega.%
\end{array}
\right. \eqno{(AS)_R}
\end{equation*}
Without loss of generality, we will consider that $\alpha=\beta=0$. Our
result main related to $(AS)_R$ is the following:

\begin{theorem}
\label{TA1} Assume that there are $(\underline{u},\overline{v}),(\overline{u}%
,\overline{u}) \in (C^{2}(\Omega) \cap L^{\infty}(\Omega))^2$ verifying: 
\begin{equation*}
\underline{u} \leq \overline{u} \,\,\, \mbox{and} \,\,\, \underline{v} \leq 
\overline{v} \,\,\, \mbox{in} \,\,\ \overline{\Omega},
\end{equation*}
\begin{equation*}
\underline{u} \leq \alpha, \underline{v} \leq \beta \,\,\, \mbox{and} \,\,\, 
\overline{u} \geq \alpha , \overline{v} \geq \beta \,\,\, \mbox{on} \,\,\
\partial \Omega
\end{equation*}
and that , 
\begin{equation*}
\left\{ 
\begin{array}{l}
\Delta \underline{u}+ b_1(x)\xi_R(|\nabla \underline{u} |^{q_1}) \geq F_u(x,\underline{u}%
, \underline{v}) \,\,\, \mbox{in} \,\,\, \Omega, \\ 
\mbox{} \\ 
\Delta \underline{v}+ b_2(x)\xi_R(|\nabla \underline{v}|^{q_2}) \geq F_v(x,\underline{u}%
, \underline{v}) \,\,\, \mbox{in} \,\,\, \Omega,%
\end{array}
\right.
\end{equation*}
and 
\begin{equation*}
\left\{ 
\begin{array}{l}
\Delta \overline{u}+ b_1(x)\xi_R(|\nabla \overline{u}|^{q_1}) \leq F_u(x,\overline{u}, 
\overline{v}), \,\,\, \mbox{in} \,\,\, \Omega, \\ 
\mbox{} \\ 
\Delta \underline{v}+ b_2(x)\xi_R(|\nabla \overline{v}|^{q_2}) \leq F_v(x,\overline{u}, 
\overline{v}), \,\,\, \mbox{in} \,\,\, \Omega.%
\end{array}
\right.
\end{equation*}
Then, there is $(u,v) \in (H^{1}(\Omega))^2$ such that 
\begin{equation*}
\underline{u} \leq u \leq \overline{u} \,\,\, \mbox{and} \,\,\, \underline{v}
\leq v \leq \overline{v} \,\,\, \mbox{in} \,\,\ \overline{\Omega}
\end{equation*}
and 
\begin{equation*}
\left\{ 
\begin{array}{l}
\Delta u+ b_1(x)\xi_R(|\nabla u|^{q_1}) = F_u(x,u,v) \,\,\, \mbox{in} \,\,\,
\Omega, \\ 
\mbox{} \\ 
\Delta v+ b_2(x)\xi_R(|\nabla v|^{q_2}) = F_v(x,u,v) \,\,\, \mbox{in} \,\,\,
\Omega, \\ 
\mbox{} \\ 
u= \alpha, v=\beta \,\,\, \mbox{on} \,\,\, \partial \Omega.%
\end{array}
\right.
\end{equation*}
\end{theorem}

\noindent \textbf{Proof.} Here, we will use a result due to Alves
\& Moussaoui \cite[Theorem 2.1]{AlvesAbdelkrim}. First of all, we observe
that without loss of generality we can assume that $\alpha=\beta=0$.
Setting the functions $H,G: \Omega \times \mathbb{R}^{+} \times \mathbb{R}^{+} \times  \mathbb{R}^{N} \times \mathbb{R}^{N} \to \mathbb{R}$ given by 
$$
H(x,s,t,\eta,\zeta)= -F_{s}(x,s,t)+b_1(x)\xi_R(|\eta|)
$$ 
and
$$
G(x,s,t,\eta,\zeta)= -F_{t}(x,s,t)+b_2(x)\xi_R(|\zeta|),
$$
we observe that they are continuous functions and given $T,S>0$, there exists $C=C(R)>0$ such that
$$
|H(x,s,t,\eta,\zeta)|,|G(x,s,t,\eta,\zeta)| \leq C, \,\,\, \forall (x,s,t,\eta,\xi) \in \Omega \times [0,T] \times [0,S] \times  \mathbb{R}^{N} \times \mathbb{R}^{N},
$$
finishing the proof of Theorem \ref{TA1}. \hfill \rule{2mm}{2mm}

\begin{corollary}
\label{C1} On the hypotheses of Theorem \ref{TA1}, if the $L^{\infty
}(\Omega )$ norms of the pairs $(\underline{u},\underline{v})$, $(\overline{u%
},\overline{v})$ are independent of $R$, for $R$ large enough, then there is 
$R^{\ast }>0$ such that if $R>R^{\ast }$, the solution $(u,v)$ given by
Theorem \ref{TA1} verifies 
\begin{equation*}
\max_{x\in \overline{\Omega }}|\nabla u(x)|,\max_{x\in \overline{\Omega }%
}|\nabla v(x)|<\min \{R^{\frac{1}{q_{1}}},R^{\frac{1}{q_{2}}}\}.
\end{equation*}%
Thus, $(u,v)$ is a solution of the original system $(S)$.
\end{corollary}

\noindent \textbf{Proof.} A first point that we should mention is the fact
that by Elliptic Regularity, 
\begin{equation*}
u,v\in W^{2,p}(\Omega )\,\,\,\forall p\in \lbrack 1,+\infty ),
\end{equation*}%
because $\xi _{R}\in L^{\infty }([0,+\infty ))$, $u$, $v\in L^{\infty
}(\Omega )$ and $F_{u},F_{v}$ are continuous functions. From now on, we will
fix $p$ such that 
\begin{equation}
W^{2,p}(\Omega )\hookrightarrow C^{1,\alpha }(\overline{\Omega })
\label{CC1}
\end{equation}%
is a continuous embedding. Now, we observe that $u$ is a solution of the
problem 
\begin{equation*}
\Delta u-u=B_{R}(x)(1+|\nabla u|^{2}),
\end{equation*}%
where 
\begin{equation*}
B_{R}(x)=\frac{-u+F_{u}(x,u,v)-b_1(x)\xi _{R}(|\nabla u|^{q_{1}})}{1+|\nabla
u|^{2} }.
\end{equation*}%
Once that 
\begin{equation*}
\xi _{R}(t)\leq t\,\,\,\forall t\geq 0,
\end{equation*}%
a direct computation shows that there is $C^{\ast }>0$, independent of $R>R^*$,
such that 
\begin{equation*}
|B_{R}(x)|\leq C^{\ast }\,\,\,\forall x\in \Omega,
\end{equation*}%
leading to 
\begin{equation}
\Vert B_{R}\Vert _{\infty }\leq C^{\ast }\,\,\,\forall R>R^*.  \label{CC2}
\end{equation}%
By using a result due to Amann \& Crandall \cite[Lemma 4]{Amann}, there is
an increasing function $\gamma _{0}:[0,+\infty )\rightarrow \lbrack 0,\infty
)$, depending only of $\Omega $, $p$ and $N$, and satisfying 
\begin{equation*}
\Vert u\Vert _{W^{2,p}(\Omega )}\leq \gamma _{0}(\Vert B_{R}\Vert _{\infty
}).
\end{equation*}%
Combining the last inequality with (\ref{CC1}) and (\ref{CC2}), 
\begin{equation*}
\Vert u\Vert _{C^{1,\alpha }(\overline{\Omega} )}\leq C\gamma _{0}(C^{\ast
}),
\end{equation*}%
for some $C>0$. Fixing 
\begin{equation*}
K=C\gamma _{0}(C^{\ast }),
\end{equation*}%
we derive that 
\begin{equation*}
|\frac{\partial u(x)}{\partial x_{i}}|\leq K\,\,\,\forall x\in \overline{%
\Omega }\,\,\ \mbox{and}\,\,\,i=1,2,..,N.
\end{equation*}%
Thereby, 
\begin{equation*}
|\nabla u(x)|\leq NK\,\,\,\forall x\in \overline{\Omega },
\end{equation*}%
implying that 
\begin{equation*}
\max_{x\in \overline{\Omega }}|\nabla u(x)|\leq NK.
\end{equation*}%
Fixing $R_{1}^{\ast }=(NK)^{q_{1}}$, it follows that 
\begin{equation*}
\max_{x\in \overline{\Omega }}|\nabla u(x)|\leq (R_{1}^{\ast })^{\frac{1}{%
q_{1}}}.
\end{equation*}%
By a similar argument, we get $R_{2}^{\ast }>0$ verifying 
\begin{equation*}
\max_{x\in \overline{\Omega }}|\nabla v(x)|\leq (R_{2}^{\ast })^{\frac{1}{%
q_{2}}}.
\end{equation*}%
Now, the corollary follows setting $R^{\ast }=\max \{R_{1}^{\ast
},R_{2}^{\ast }\}$. \hfill \rule{2mm}{2mm}

\section{Proof of Theorem \protect\ref{1.1}}

We begin the proof of Theorem \ref{1.1} by Finite Case.

\vspace{0.5 cm}

\noindent \textbf{Case 1: Finite case}

\vspace{0.5 cm}

In what follows, we fix $M>\max \left\{ \alpha ,\beta \right\} $, $m<\min
\left\{ \alpha ,\beta \right\} $ and denote by $\psi \in C^{2}\left( \Omega
\right) \cap C^{1,\alpha }(\overline{\Omega })$ the unique positive solution
of the problem 
\begin{equation*}
\left\{ 
\begin{array}{l}
\Delta \psi =\displaystyle\sum\limits_{i=1}^{2}a_{i}^{2}\left( x\right)
g\left( \psi \right) \text{ in }\Omega ,\text{ } \\ 
\psi >0\text{ in }\Omega , \\ 
\psi =m\text{ on }\partial \Omega \text{, }%
\end{array}%
\right.
\end{equation*}%
which exists by a result found in \cite[Proposition 1]{lair}. The pairs 
\begin{equation*}
(\underline{u},\underline{v})=(\psi (x),\psi (x))\,\,\,\mbox{and}\,\,\,(%
\overline{u},\overline{v})=(M,M),
\end{equation*}%
satisfy the hypotheses of Corollary \ref{C1}. Thus, the system 
\begin{equation*}
\left\{ 
\begin{array}{l}
\Delta u+b_1(x)|\nabla u|^{q_{1}}=F_{u}(x,u,v),\,\,\,\mbox{in}\,\,\,\Omega
\\ 
\mbox{} \\ 
\Delta v+b_2(x)|\nabla v|^{q_{2}}=F_{v}(x,u,v),\,\,\,\mbox{in}\,\,\,\Omega
\\ 
\mbox{} \\ 
u=\alpha ,v=\beta \,\,\,\mbox{on}\,\,\,\partial \Omega ,%
\end{array}%
\right. \eqno{(S)_{\alpha,\beta}}
\end{equation*}%
has a solution. Moreover, by elliptic regularity, we must have $u,v\in
C^{2}(\Omega )\cap C^{1,\alpha }(\overline{\Omega })$. \hfill\rule{2mm}{2mm}

\vspace{0.5 cm}

\noindent \textbf{Case 2: Infinite case}.

\vspace{0.5 cm}

In this case, we denote by $(u_{n},v_{n})$ the solution of the system 
\begin{equation}
\left\{ 
\begin{array}{l}
\Delta u+b_{1}\left( x\right) \left\vert \nabla u\right\vert
^{q_{1}}=F_{u}\left( x,u,v\right) \text{ in }\Omega , \\ 
\mbox{} \\ 
\Delta v+b_{2}\left( x\right) \left\vert \nabla v\right\vert
^{q_{2}}=F_{v}\left( x,u,v\right) \text{ in }\Omega , \\ 
\mbox{} \\ 
u=v=n\text{ on }\partial \Omega ,%
\end{array}%
\right.  \label{in}
\end{equation}%
which exists by finite case. We remark that $(u_{n},v_{n})$ can be chosen
satisfying the inequality%
\begin{equation}
u_{n}\leq u_{n+1}\,\text{and }v_{n}\leq v_{n+1}\,\,\forall n\in \mathbb{N}.
\label{3.1}
\end{equation}%
Indeed, note that $(u_{1},v_{1})$ satisfies 
\begin{equation}
\left\{ 
\begin{array}{l}
\Delta u+b_{1}\left( x\right) \xi _{R}(\left\vert \nabla u\right\vert
^{q_{1}})=F_{u}\left( x,u,v\right) \text{ in }\Omega , \\ 
\mbox{} \\ 
\Delta v+b_{2}\left( x\right) \xi _{R}(\left\vert \nabla v\right\vert
^{q_{2}})=F_{v}\left( x,u,v\right) \text{ in }\Omega , \\ 
\mbox{} \\ 
u_{1}=v_{1}=1\leq 2\text{ on }\partial \Omega ,%
\end{array}%
\right.  \label{i2}
\end{equation}%
for all 
\begin{equation*}
R>\max \{\max_{x\in \overline{\Omega }}|\nabla u_{1}(x)|^{q_{1}},\max_{x\in 
\overline{\Omega }}|\nabla v_{1}(x)|^{q_{2}}\},
\end{equation*}%
because 
\begin{equation*}
\xi _{R}(\left\vert \nabla u_1\right\vert ^{q_{1}}(x))=\left\vert \nabla u_1 (x)
\right\vert ^{q_{1}} \,\,\,\mbox{and}\,\,\,\xi _{R}(\left\vert \nabla
v_1(x)\right\vert ^{q_{2}})=\left\vert \nabla v_1(x)\right\vert
^{q_{2}}\,\,\,\forall x\in \overline{\Omega }.
\end{equation*}%
Applying Corollary \ref{C1} with $(\underline{u},\underline{v}%
)=(u_{1},v_{1}) $ and $(\overline{u},\overline{v})=(2,2)$ , there exists a
solution $(u_{2},v_{2})$ of 
\begin{equation*}
\left\{ 
\begin{array}{l}
\Delta u+b_{1}(x)\left\vert \nabla u\right\vert ^{q_{1}}=F_{u}(x,u,v)\text{
in }\Omega , \\ 
\mbox{} \\ 
\Delta v+b_{2}(x)\left\vert \nabla v\right\vert ^{q_{2}}=F_{v}(x,u,v)\text{
in }\Omega , \\ 
\mbox{} \\ 
u_{1}=v_{2}=2\text{ on }\partial \Omega ,%
\end{array}%
\right.
\end{equation*}%
satisfying 
\begin{equation*}
u_{1}\leq u_{2}\,\,\,\mbox{and}\,\,\ v_{1}\leq v_{2}\leq M_{1}:=2\,\,\,\text{
in }\,\,\,\overline{\Omega }.
\end{equation*}%
Repeating the above argument, of an iterative way, for each $M_{n}=n+1$; $%
n=1,2,...$, the pair $\left( u_{n},v_{n}\right) $ satisfies 
\begin{equation*}
\left\{ 
\begin{array}{l}
\Delta u_{n}+b_{1}\left( x\right) \left\vert \nabla u_{n}\right\vert
^{q_{1}}=F_{u}\left( x,u_{n},v_{n}\right) \text{ in }\Omega , \\ 
\mbox{} \\ 
\Delta v+b_{2}\left( x\right) \left\vert \nabla v\right\vert
^{q_{2}}=F_{v}\left( x,u_{n},v_{v}\right) \text{ in }\Omega , \\ 
\mbox{} \\ 
u_{n}=v_{n}\leq n+1\text{ on }\partial \Omega .%
\end{array}%
\right.
\end{equation*}%
Applying again Corollary \ref{C1} with $(\underline{u},\underline{v}%
)=(u_{n},v_{n})$ and $(\overline{u},\overline{v})=(n+1,n+1)$, we get a
solution $(u_{n+1},v_{n+1})$ of 
\begin{equation*}
\left\{ 
\begin{array}{l}
\Delta u+b_{1}\left( x\right) \left\vert \nabla u\right\vert
^{q_{1}}=F_{u}\left( x,u,v\right) \text{ in }\Omega , \\ 
\mbox{} \\ 
\Delta v+b_{2}\left( x\right) \left\vert \nabla v\right\vert
^{q_{2}}=F_{v}\left( x,u,v\right) \text{ in }\Omega , \\ 
\mbox{} \\ 
u=v=n+1\text{ on }\partial \Omega ,%
\end{array}%
\right.
\end{equation*}%
with 
\begin{equation*}
u_{n}\leq u_{n+1}\,\,\,\mbox{and}\,\,\ v_{n}\leq v_{n+1}\leq M_{n}:=n+1\,\,\,%
\text{ in }\,\,\,\overline{\Omega }.
\end{equation*}%
Once that sequences $(u_{n})$ and $(v_{n})$ are nondecreasing, there are
functions $u$, $v:\Omega \rightarrow \mathbb{R}$ verifying 
\begin{equation*}
u_{n}(x)\rightarrow u(x)\,\,\,\mbox{and}\,\,\,v_{n}(x)\rightarrow v(x)\,\,\,%
\mbox{in}\,\,\,\Omega .
\end{equation*}%
From now on, we denote by $\widetilde{u}$ and $\widetilde{v}$ the solutions
of the problems 
\begin{equation*}
\left\{ 
\begin{array}{l}
\Delta u+\left\Vert b_{1}\right\Vert _{{\infty }}\left\vert \nabla
u\right\vert ^{q_{1}}=\underset{x\in \overline{\Omega }}{\min }a_{1}\left(
x\right) f_{1}\left( u\right) \text{ in }\Omega, \\ 
u>0\text{ in }\Omega, \\ 
u=+\infty \text{ on }\,\,\partial \Omega,%
\end{array}%
\right.
\end{equation*}%
and 
\begin{equation*}
\left\{ 
\begin{array}{l}
\Delta v+\left\Vert b_{2}\right\Vert _{{\infty }}\left\vert \nabla
v\right\vert ^{q_{2}}=\underset{x\in \overline{\Omega }}{\min }a_{2}\left(
x\right) f_{2}\left( v\right) \text{ in }\Omega, \\ 
v>0\text{ in }\Omega, \\ 
v=+\infty \,\,\,\mbox{on}\,\,\,\partial \Omega ,%
\end{array}%
\right.
\end{equation*}%
which exist from a result due to Bandle \& Giarrusso \cite{bg}. We claim
that 
\begin{equation}
(i)\,\,u_{n}\leq \widetilde{u}\,\,\,\mbox{and}\,\,\,(ii)\,\,v_{n}\leq 
\widetilde{v}\,\,\,\text{ in }\,\,\,\overline{\Omega }\,\,\,\forall n\geq 1.
\label{F1}
\end{equation}%
Indeed, suppose by contradiction that (\ref{F1})(i) does not hold. Then,
there exists $x_{0}\in \Omega $ such that 
\begin{equation*}
u_{n}\left( x_{0}\right) >\widetilde{u}\left( x_{0}\right) \text{ in }%
\,\,\Omega \,\,\text{ for some }\,\,n\geq 1.
\end{equation*}%
Since 
\begin{equation*}
\lim_{d(x, \partial \Omega) \to 0 }\left[ u_{n}\left( x\right) -\widetilde{u}%
\left( x\right) \right] =-\infty ,
\end{equation*}%
we deduce that $\displaystyle \max_{x\in \overline{\Omega }}(u_{n}-%
\widetilde{u})(x) $ is achieved in $\Omega $, for example at $x_{1}\in
\Omega $. This form, 
\begin{equation*}
\left\vert \nabla u_{n}\left( x_{1}\right) \right\vert ^{q_{1}}=\left\vert
\nabla \widetilde{u}\left( x_{1}\right) \right\vert ^{q_{1}}
\end{equation*}%
and 
\begin{eqnarray*}
0 &\geq &\Delta \left( u_{n}-\widetilde{u}\right)(x_1) \\
&=&a_{1}\left( x_{1}\right) f_{1}\left( u_{n}\right) -b_{1}\left(
x_{1}\right) \left\vert \nabla u_{n}\right\vert ^{q_{1}}-\underset{x\in 
\overline{\Omega} }{\min }a_{1}\left( x_{1}\right) f_{1}\left( \widetilde{u}%
\right) +\left\Vert b_{1}\right\Vert _{{\infty }}\left\vert \nabla 
\widetilde{u}\right\vert ^{q_{1}}>0,
\end{eqnarray*}%
obtaining a contradiction. Therefore, (\ref{F1})(i) holds. The same argument
works to prove (\ref{F1})(ii).

On the other hand, by using a well known result due to Ladyzenskaya and
Ural'treva \cite{olga}, given $\Omega _{1}\subset \subset \Omega _{2}\subset
\subset \Omega $, there is $C>0$ such that 
\begin{equation}
\max_{x\in \overline{\Omega _{1}}}|\nabla u_{n}|\leq C\max_{x\in \overline{%
\Omega _{2}}}|u_{n}|\leq C\max_{x\in \overline{\Omega _{2}}}|\widetilde{u}%
|=K_{5}\,\,\ \forall n\in \mathbb{N}  \label{F2}
\end{equation}%
and 
\begin{equation}
\max_{x\in \overline{\Omega _{1}}}|\nabla v_{n}|\leq C\max_{x\in \overline{%
\Omega _{2}}}|v_{n}|\leq C\max_{x\in \overline{\Omega _{2}}}|\widetilde{v}%
|=K_{6}\,\,\ \forall n\in \mathbb{N}.  \label{F3}
\end{equation}%
Combining (\ref{F2}), (\ref{F3}) and Elliptic Regularity, it follows that
there are subsequences of $(u_{n})$ and $(v_{n})$, still denoted by $(u_{n})$
and $(v_{n})$, such that 
\begin{equation*}
u_{n}\rightarrow u\,\,\ \mbox{and}\,\,\ v_{n}\rightarrow v\,\,\ \mbox{in}%
\,\,\ C_{loc}^{2}(\Omega ).
\end{equation*}%
This fact yields $u,v\in C^{2}(\Omega )$ and 
\begin{equation*}
\left\{ 
\begin{array}{l}
\Delta u+b_{1}\left( x\right) \left\vert \nabla u\right\vert
^{q_{1}}=F_{u}\left( x,u,v\right) \text{ in }\Omega , \\ 
\mbox{} \\ 
\Delta v+b_{2}\left( x\right) \left\vert \nabla v\right\vert
^{q_{2}}=F_{v}\left( x,u,v\right) \text{ in }\Omega , \\ 
\mbox{} \\ 
u,v>0\,\,\text{ in }\,\,\Omega .%
\end{array}%
\right.
\end{equation*}%
To complete the proof, it suffices to prove that $(u,v)$ blows up at the
boundary. Arguing by contradiction, we will assume that $u$ does not blow up
at the boundary. Then, there exist $x_{0}\in \partial \Omega $ and $\left(
x_{k}\right) \subset \Omega $\ such that 
\begin{equation*}
\lim_{k\rightarrow \infty }x_{k}=x_{0}\text{ and }\lim_{k\rightarrow \infty
}u\left( x_{k}\right) =L\in \left( 0,\infty \right) .
\end{equation*}%
In what follows, fix $n>4L$ and $\delta >0$ such that $u_{n}\left( x\right)
\geq n/2$ for all $x\in \overline{\Omega }_{\delta },$ where 
\begin{equation*}
\overline{\Omega }_{\delta }=\left\{ x\in \overline{\Omega }\left\vert
dist\left( x,\partial \Omega \right) \leq \delta \right. \right\} .
\end{equation*}%
Then, for $k$ large enough, $x_{k}\in \overline{\Omega }_{\delta }$ and $%
u_{n}\left( x_{k}\right) >2L$. Since 
\begin{equation*}
u_{n}\left( x_{k}\right) \leq u_{n+1}\left( x_{k}\right) \leq ...\leq
u_{n+j}\left( x_{k}\right) \leq ...\leq u\left( x_{k}\right) \forall j\text{%
, }
\end{equation*}%
we have that $u\left( x_{k}\right) \geq 2L,$ which is a contradiction.
Therefore, $u$ blows up at the boundary. The same approach can be used to
prove that $v$ also blows up. \hfill \rule{2mm}{2mm}

\vspace{0.5 cm}

\noindent \textbf{Case 3: Semifinite case}

\vspace{0.5 cm}

Let $(u_{n},v_{n})\in (C^{2}(\Omega )\cap C^{1,\alpha }(\overline{\Omega }%
))^{2}$ be a solution of $(S)_{\alpha ,\beta }$ with $\alpha =n$, $n\in 
\mathbb{N}$ and $\beta $ fixed. As in the previous case, the sequence $%
(u_{n},v_{n})$ is bounded on compact subset contained in $\Omega $, implying
that there exist functions $u$, $v$ verifying 
\begin{equation*}
u_{n}\rightarrow u~~\mbox{in}~~C^{2}(K)
\end{equation*}%
and 
\begin{equation*}
v_{n}\rightarrow v~~\mbox{in}~~C^{2}(K)
\end{equation*}%
for any compact subset $K\subset \Omega $. Moreover, the arguments used in
the previous cases give that $u$ blows up at the boundary, that is, 
\begin{equation*}
u(x)=+\infty ~~~\mbox{on}~~\partial \Omega .
\end{equation*}%
Related to sequence $(v_{n})$, we recall that 
\begin{equation*}
\left\{ 
\begin{array}{lll}
\Delta v_{n}+b_{2}(x)|\nabla v_{n}|^{q_{2}}=F_{v}(x,u_{n},v_{n})\,\,\,~~%
\mbox{in}~~\Omega , &  &  \\ 
\mbox{} &  &  \\ 
v_{n}>0\,\,\,\mbox{in}\,\,\,\Omega, &  &  \\ 
\mbox{} &  &  \\ 
v_{n}=\beta \,\,\,\mbox{on}~~\partial \Omega, &  & 
\end{array}%
\right.
\end{equation*}%
\noindent and 
\begin{equation*}
v_{n}(x)\leq \beta ~~\mbox{in}~~\forall x\in \overline{\Omega }~~\mbox{and}%
~~n\geq 1.
\end{equation*}%
\noindent Passing to the limit as $n\rightarrow +\infty $, we obtain that $%
v(x)\leq \beta $ for all $x\in \Omega $.

\begin{claim}
\label{A3} Let $x_0 \in \partial \Omega$ and $(x_k) \subset \Omega$ be a
sequence with $x_{k} \to x_0$. Then $v(x_k) \to \beta$ as $k \to +\infty$.
\end{claim}

Indeed, if the limit does not hold, there exist $\epsilon >0$ and a
subsequence of $(x_k)$, still denoted by itself, such that 
\begin{equation}  \label{X1}
x_k\rightarrow x_0 ~~ \mbox{and} ~~ v(x_k) \leq \beta-\varepsilon ~~ \forall
k \in \mathbb{N}.
\end{equation}
Since $v_1=\beta$ on $\partial\Omega$, there is $\delta >0$ such that 
\begin{equation*}
v_1(x)\geq\beta-\frac{\varepsilon}{2}, ~~ \forall x \in \overline{\Omega}%
_{\delta}.
\end{equation*}
Hence, for $k$ large enough, 
\begin{equation*}
x_k \in \overline{\Omega}_{\delta}
\end{equation*}
and 
\begin{equation*}
v_1(x_k)\geq\beta-\frac{\varepsilon}{2}>\beta-\varepsilon.
\end{equation*}
Recalling that $v_1 \leq v$ in $\Omega$, it follows that 
\begin{equation*}
v(x_k)\geq \beta-\frac{\varepsilon}{2}>\beta-\varepsilon,
\end{equation*}
obtaining a contradiction with \eqref{X1}.

From Claim \ref{A3}, we can continuously extend the function $v$ from $%
\Omega $ to $\overline\Omega$, by considering 
\begin{equation*}
v(x)=\beta ~~ \mbox{on} ~~ \partial\Omega,
\end{equation*}
concluding this way the proof of the semifinite case. \hfill \rule{2mm}{2mm}

\subsection{Proof of Theorem \protect\ref{entire}}

Firstly, we provide a lower bound for the system $(LS)$. To this end, we
consider the function $w:\mathbb{R}^{N}\rightarrow \left[ 0,\infty \right) $
implicitly defined by%
\begin{equation*}
z\left( x\right) =\int_{w\left( x\right) }^{\infty }\frac{1}{g\left(
t\right) }dt \,\,\,\,\, \text{ }x\in \mathbb{R}^{N},
\end{equation*}%
where $z$ was given in (\ref{1.6.}). This is possible due to the fact that
Keller-Osserman (\ref{ko}) condition gives 
\begin{equation*}
\int_{1}^{\infty }\frac{1}{g\left( t\right) }dt<\infty. \,\,\,~~~~~~~%
\mbox{(see \cite{covei}, \cite{lair} for details).}
\end{equation*}%
Note that $w\in C^{2}\left( \mathbb{R}^{N},\left( 0,\infty \right) \right) $%
, $w\left( x\right) \rightarrow +\infty $ as $\left\vert x\right\vert
\rightarrow \infty $ and%
\begin{equation*}
\Delta w\left( x\right) \geq \sum_{i=1}^{2}a_{i}^{2}\left( x\right) g\left(
w\left( x\right) \right) \text{ for all }x\in \mathbb{R}^{N}.
\end{equation*}%
Moreover, 
\begin{equation*}
\Delta w\geq F_{u}\left( x,w,w\right) \text{ in }B_{n}\newline
\text{, }w\leq w_{n}\text{ on }\partial B_{n},\newline
\end{equation*}%
and 
\begin{equation*}
\Delta w\geq F_{v}\left( x,w,w\right) \text{ in }B_{n}\newline
\text{, }w\leq w_{n}\text{ on }\partial B_{n}.\newline
\end{equation*}%
Using function $w$, we consider the system 
\begin{equation}
\left\{ 
\begin{array}{l}
\Delta u+b_{1}\left( x\right) \left\vert \nabla u\right\vert
^{q_{1}}=F_{u}\left( x,u,v\right) \,\, \text{ in } \,\, B_{n}, \\ 
\mbox{} \\ 
\Delta v+b_{2}\left( x\right) \left\vert \nabla v\right\vert
^{q_{2}}=F_{v}\left( x,u,v\right) \,\, \text{ in } \,\, B_{n}, \\ 
\mbox{} \\ 
u,v>0\,\,\,\mbox{in} \,\, \Omega, \\ 
\mbox{} \\ 
u=v=w_{n} \,\, \text{ on } \,\, \partial B_{n},%
\end{array}%
\right.  \label{ball}
\end{equation}%
where $B_{n}$ is the open ball of radius $n$ centered at the origin and 
\begin{equation*}
w_{n}=\max_{x\in \overline{B}_{n}}w\left( x\right) .
\end{equation*}%
Applying Theorem \ref{1.1} with $(\underline{u},\underline{u})=(w,w)$ and $(%
\overline{u},\overline{v})=(w_{n},w_{n})$, there is a solution $%
(u_{n},v_{n})\in \left[ C^{2}(B_{n}\right) \cap C^{1,\alpha }(\overline{%
\Omega })]^{2}$ of (\ref{ball}). Moreover, we can choose the sequence $%
(u_n,v_n)$ satisfying 
\begin{equation*}
w(x)\leq u_{n}(x)\leq u_{n+1}(x)\,\,\,\text{ forall }x\in \overline{B}_{n},
\end{equation*}%
and 
\begin{equation*}
w(x)\leq v_{n}(x)\leq v_{n+1}(x)\,\,\,\text{ forall }x\in \overline{B}_{n}.
\end{equation*}%
From this, there are $u,v:\mathbb{R}^{N}\rightarrow \mathbb{R}$ such that 
\begin{equation*}
u_{n}(x)\rightarrow u(x)\,\,\,\mbox{and}\,\,\,v_{n}(x)\rightarrow
v(x)\,\,\,\forall x\in \mathbb{R}^{N}.
\end{equation*}%
Arguing as in the previous sections, there are subsequences of $%
(u_{n}),(v_{n})$, still denoted by themself, such that 
\begin{equation*}
u_{n}\rightarrow u\,\,\ \mbox{and}\,\,\,v_{n}\rightarrow v\,\,\,\mbox{in}%
\,\,\,C_{loc}^{2}(\mathbb{R}^{N})
\end{equation*}%
and 
\begin{equation*}
u(x),v(x)\geq w(x)\,\,\,\forall x\in \mathbb{R}^{N}.
\end{equation*}%
Consequently, $u,v\in C^{2}\left( \mathbb{R}^{N}\right) $ and $(u,v)$ is a
solution of 
\begin{equation*}
\left\{ 
\begin{array}{l}
\Delta u+b_{1}\left( x\right) \left\vert \nabla u\right\vert
^{q_{1}}=F_{u}(x,u,v)\text{ in }\mathbb{R}^{N}, \\ 
\mbox{} \\ 
\Delta v+b_{2}(x)\left\vert \nabla v\right\vert ^{q_{2}}=F_{v}(x,u,v)\text{
in }\mathbb{R}^{N}, \\ 
\mbox{} \\ 
u,v>0\,\,\text{ in }\mathbb{R}^{N}, \\ 
\mbox{} \\ 
u(x),v(x)\rightarrow +\infty \,\,\,\mbox{as}\,\,\,|x|\rightarrow +\infty ,%
\end{array}%
\right.
\end{equation*}%
showing that $(u,v)$ is a entire large solution for $(LS)$. \hfill\rule%
{2mm}{2mm}

\end{document}